\documentclass{amsart}
\usepackage{amssymb}
\input xy
\xyoption{all}
\newcommand{\ra}{\rightarrow}
\newcommand{\str}{\mathcal{O}}

\newcommand{\Ps}{\mathbf{P}}
\newcommand{\Q}{\mathbf{Q}}
\newcommand{\F}{\mathbf{F}}
\newcommand{\C}{\mathbf{C}}
\newcommand{\Z}{\mathbf{Z}}
\newcommand{\cC}{\mathcal{C}}
\newcommand{\cU}{\mathcal{U}}
\newcommand{\cF}{\mathcal{F}}
\newcommand{\cV}{\mathcal{V}}

\renewcommand{\phi}{\varphi}
    \newtheorem{Lem}{Lemma}[section]
    \newtheorem{Prop}[Lem]{Proposition}
    \newtheorem{Thm}[Lem]{Theorem}
    
    \newtheorem{Cor}[Lem]{Corollary}

   \theoremstyle{definition}
    \newtheorem{Def}[Lem]{Definition}

    \newtheorem{Rem}[Lem]{Remark}
    \newtheorem{Dia}[Lem]{Diagram}
\DeclareMathOperator{\Fr}{Fr}
\DeclareMathOperator{\Tr}{Tr}
\DeclareMathOperator{\tor}{tor}

\DeclareMathOperator{\Sym}{Sym}
\DeclareMathOperator{\et}{\acute{e}t}

\DeclareMathOperator{\PGL}{PGL}
\DeclareMathOperator{\GL}{GL}
\DeclareMathOperator{\SL}{SL}
\DeclareMathOperator{\Gal}{Gal}
\DeclareMathOperator{\Gr}{Gr}
\DeclareMathOperator{\spec}{Spec}
\DeclareMathOperator{\rank}{rank}
\DeclareMathOperator{\Pic}{Pic}
\DeclareMathOperator{\Aut}{Aut}
\DeclareMathOperator{\red}{red}

\begin{document}
\title{The $L$-series of a cubic fourfold}
\author{Klaus Hulek}
\address{Institut f\"ur Algebraische Geometrie, Universit\"at Hannover, Welfengarten 1, D-30167, Hannover, Germany}
\email{hulek@math.uni-hannover.de}

\author{Remke Kloosterman}
\address{Institut f\"ur Algebraische Geometrie, Universit\"at Hannover, Welfengarten 1, D-30167, Hannover, Germany}
\email{kloosterman@math.uni-hannover.de}

\begin{abstract}
{We study the $L$-series of cubic fourfolds. Our main result is that, if $X/\C$ is a special 
cubic fourfold associated to some polarized $K3$ surface $S$, defined over a number field $K$ and satisfying $S^{[2]}(K)\neq \emptyset$, then $X$ has a model over $K$ such that the $L$-series of the primitive cohomology of $X/K$ can be expressed in the $L$-series of $S/K$.
This allows us to compute the $L$-series for a discrete dense
subset of cubic fourfolds in the moduli spaces of certain special cubic fourfolds. We also discuss a concrete
example.
}
\end{abstract}
\thanks{This work was partially supported by the DFG Schwerpunktprogramm ``Globale Methoden
in der komplexen Geometrie'' under grant HU 337/5-3. It was started while the first 
named author enjoyed the hospitality of the Fields Institute in Toronto. 
We are grateful to a number of mathematicians whose
help was crucial to us, in particular to V.~Batyrev whose question about the $L$-series of cubic fourfolds 
started this work and to 
B.~Hassett who taught us about about the geometry of cubic fourfolds, in particular in 
connection with our example. We also profited from discussions or exchanges of e-mails with
M.~Artebani,
A.~Beauville, S.~Kondo, A.~Laface, M.~Lehn, S.~Mukai, K.~Oguiso,  S.~Popsecu, A.~Sarti, M.~Sch\"utt and D.-Q.~Zhang. 
}
\maketitle

\section{Introduction}\label{intro}

After the proof of the Taniyama-Shimura-Weil conjecture by Wiles, Taylor and others, many efforts were made to prove similar results for
other classes of varieties. The most natural class to investigate is that of Calabi-Yau
varieties. In the case of $K3$ surfaces Shioda and Inose \cite{SI} had already determined the zeta-function of
singular $K3$ surfaces defined over some number field $K$. In the case that the singular $K3$ surface
$S$ is defined over $\Q$, Livn\'e \cite{LivMot} has shown that the $L$-series associated to
the $2$-dimensional transcendental lattice of $S$ comes from a weight $3$ modular form with complex
multiplication in some quadratic imaginary field $L$.
For rigid Calabi-Yau threefolds defined over $\Q$, Dieulefait and
Manoharmayum \cite{DieMan} have proved modularity under mild restrictions on the primes of bad reduction.
Further examples and results are now also known for certain non-rigid Calabi-Yau threefolds and some
higher dimensional examples. For a survey see \cite{HS} .

The main purpose of this paper is to show that there are many cubic fourfolds whose $L$-series can be computed.
A cubic fourfold $X$ is called special if it contains a surface $T$ which is not homologous to a multiple
of the class $h^2$ where $h$ is the class of the hyerplane section. The discriminant of $X$ is the discriminant
of the saturated rank $2$ sublattice of $H^4(X,\Z)$ spanned by $h^2$ and $T$. These cubic fourfolds were
studied extensively by Hassett \cite{HasArt} . He proved the following result: if the discriminant $d$ is
greater than $6$ and if $ d \equiv 0,2 \bmod 6$,
then the special fourfolds of discriminant $d$ are parameterized by a non-empty $19$-dimensional quasi-projective
variety $\cC_d$. These cubic fourfolds are closely related to $K3$ surfaces, more precisely if $d=2(n^2+n+1)$ for some $n\in \Z_{>1}$, then there exists a
non-empty open subset of $\cC_d$ such that the Fano variety $F(X)$ of a cubic fourfold belonging to this
open set is isomorphic to the desingularized second symmetric product $S^{[2]}$ of some $K3$ surface $S$.

Fix such a $K3$-surface $S$. Consider the composition \[\psi:S^{[2]}\hookrightarrow \Gr(1,5)\hookrightarrow \Ps^{14}.\]
The first inclusion is given by $F(X)\cong S^{[2]}$, the second by the Pl\"ucker embedding.  To $\psi$ we can associate a line bundle $\mathcal{L}$. 

It turns out that $\mathcal{L}$ is isomorphic to  $\str(n\Delta+f)$, where $\Delta$ is the exceptional divisor on $S^{[2]}$ and $\str(f)$ is a line bundle on $S$ pushed forward via the diagonal embedding. We call $\str(f)$ the {\em associated line bundle}  to the isomorphism $F(X)\cong S^{[2]}$.

Our main result (Theorem~\ref{mainthm}) now says the following:
\begin{Thm} Let $K$ be a number field. Let $S/K$ be a $K3$ surface. Suppose $S^{[2]}(K)\neq \emptyset$.
Assume there is a cubic fourfold $X/\C$ such that $F(X)\cong S^{[2]}_{\C}$ and that the associated line bundle on $S$ descends to $K$. Then $X$ has a model over $K$,
\[ F(X)\cong S^{[2]} \mbox{ and } H^4_{\et}(X,\Q_{\ell})=H^2_{\et}(S,\Q_{\ell})(1)\oplus \Q_{\ell}[\Delta](1).\]
\end{Thm}
In the case where the $K3$ surface $S$ is singular
(i.e., has Picard number $20$), we can then determine the $L$-series of $X$
explicitly (see Corollary~\ref{cormodA}).
Since the singular $K3$ surfaces are everywhere dense (in the $\C$-topology) in the $19$-dimensional moduli
space $\cF_{d}$ of degree $d$ polarized $K3$ surfaces, we have
in this way found a countable number of points in a non-empty
Zariski open subset of $\cC_{2(n^2+n+1)}$ such that the corresponding fourfolds have a model over an explicitly known
number field $K$ and where the $L$-function can be computed in terms of Hecke Gr\"ossencharakters.

The organization of this paper is as follows. In Section~\ref{sec_geo} we rephrase some results on the cohomology of cubic fourfolds, their Fano varieties, and symmetric products in terms of \'etale cohomology. In Section~\ref{sec_moduli}
we recall Hasset's results from \cite{HasArt} on special cubic fourfolds.

In Section~\ref{sec_lseries} we prove that the isomorphism $F(X)\cong S^{[2]}$ descends to any field $K$ for which $S^{[2]}(K)\neq \emptyset$. The proof uses a classical description of morphisms onto a Grassmannian variety in terms of vector bundles and descent theory for vector bundles. In Section~\ref{sec_fermat} we relate the $L$-function of the Fermat cubic fourfold with the $L$-function of  a weight 3 Hecke eigenform.

\subsection{Convention:} Suppose $X/S$ is a scheme. Given a morphism $T\ra S$, we denote by $X_T$ the base
change of $X$ with respect to $T\ra S$. In the case that $T=\spec R $ we write $X_R$ instead of $X_{\spec{R}}$.

\section{Geometry of cubic fourfolds}\label{sec_geo}
In this section we recall basic facts about the geometry of cubic fourfolds.
Fix a number field $K$ and a cubic fourfold $X\subset \Ps^5_K$.

The general cubic fourfold contains only surfaces whose cycle class is a multiple of $h^2$,
where $h$ is the
hyperplane section. Our main interest, however, lies in so-called special cubic fourfolds.
\begin{Def}
A cubic fourfold $X$ is called \emph{(geometrically) special} if $X_{\overline{K}}$ contains an algebraic
surface $T/\overline{K}$ which is not
homologous to $h^2$, where $h$ is the class of the hyperplane section.

Let $X$ be a special cubic fourfold and $T \subset X$ a surface which is not homologous to a multiple of
$h^2$. We denote by $L_T$ the saturated sublattice of $H^4(X,\Z)$ spanned by the classes of $T$ and $h^2$.
The \emph{discriminant} $d(X,T)$ of the pair $(X,T)$ is defined as the discriminant of the
lattice $L_T$. We call
a cubic fourfold $X$ a \emph{special cubic fourfold of discriminant $d$} if it contains a surface $T$ such
that $d(X,T)=d$. \end{Def}

It was shown by Hassett \cite[Theorem 1.0.1]{HasArt} that $d > 6$ and $ d \equiv 0,2 \bmod 6$ is a necessary
and sufficient condition for the existence of special cubic fourfolds. (See also Proposition~\ref{prpExistence}.) Note that the discriminant of a special cubic fourfold is not uniquely
determined by $X$ itself. In general, $X$ can (and will) often be special with respect to more than
one discriminant.

It is a classical fact (cf. \cite{AK}) that $X$ is
covered by lines. We denote by $F(X)$ the \emph{Fano variety} of $X$, i.e.,
\[
F(X) = \{\ell \in \Gr(1,5); \ell \subset X \}.
\]
Let $U\subset \Gr(1,5)\times \Ps^5$ be the universal line and let $p$ and $q$ denote the projection
from $U$ onto the first and second factor respectively. Especially, since $X$ is covered by lines,
the morphisms $p$ and $q$ give a correspondence between $X$ and $F(X)$.

By Lefschetz' hyperplane theorem, we have that $H^i_{\et}(X,\Q_{\ell})\cong H^i_{\et}(\Ps^5,\Q_{\ell})$
for $i\neq 4$. {}From Bott's formula and the comparison theorem between sheaf and \'etale cohomology we obtain
that $H^4_{\et}(X,\Q_{\ell})$ is a 23-dimensional Galois module.
We proceed by proving some results on $H^4_{\et}(X,\Q_{\ell})$

For any Galois module $M$, we denote by $M(k)$ the $k$-th Tate twist of $M$.
\begin{Prop} \label{prpGaloisisomorphism}
We have that $p_*q^*: H^4_{\et}(X,\Q_{\ell})\ra H^2_{\et}(F(X),\Q_{\ell})(1)$ is an isomorphism of Galois
modules.\end{Prop}

{}\begin{proof}
It follows from \cite[Proposition 6]{BeDo} that the map between the
sheaf cohomology groups $p_*q^*:H^4(X_{\C},\C)\ra H^2(F(X)_{\C},\C)$ is an isomorphism, hence
has no kernel and cokernel.
{}From the comparison theorem between sheaf and \'etale cohomology it follows that $p_*q^*: H^4_{\et}(X,\Q_{\ell})\ra H^2(F(X),\Q_{\ell})$ is an isomorphism of $\Q_{\ell}$-vector spaces.

Note that $p_*$ is given by taking
the inverse of the Poincar\'e dual of the map $p^*: H^6_{\et}(F(X),\Q_{\ell}) \to H^6_{\et}(U,\Q_{\ell})$.
By Poincar\'e duality $H^4_{\et}(U,\Q_{\ell}) \cong H^6_{\et}(U,\Q_{\ell})^{\vee}(5)$ and
$H^6_{\et}(F(X),\Q_{\ell})^{\vee}(5) \cong H^2_{\et}(F(X),\Q_{\ell})(1)$ which shows that $p_*q^*: H^4_{\et}(X,\Q_{\ell})\ra H^2_{\et}(F(X),\Q_{\ell})(1)$  is Galois-equivariant. Combining this with the fact that it is an isomorphism of vector spaces, we obtain that $p_*q^*$ is an isomorphism of Galois-modules.
\end{proof}

\begin{Rem} The map $p_*q^*$ is called the {\em Abel-Jacobi map}.\end{Rem}

\begin{Rem}
Beauville and Donagi \cite{BeDo} also proved that the Fano variety $F(X)_{\C}$ is a symplectic $4$-manifold.
More precisely,
they show that $F(X)_{\C}$ is a deformation of a variety $S^{[2]}$ where $S$ is a complex $K3$ surface and
$S^{[2]}$ is the Hilbert scheme of $0$-cycles of length $2$. It is easy to see that $S^{[2]}$ is the blow-up
of the second symmetric product of the $K3$ surface $S$ along its (singular) diagonal.
Alternatively, we can blow up the product $S \times S$ along the diagonal and then take the quotient
with respect to the involution $\iota$ given by interchanging the two factors.
\end{Rem}

The above remark shows that there is a close relation between symmetric products of some $K3$ surfaces $S$ and
cubic fourfolds $X$.

Fix a $K3$-surface $S/K$ and consider the diagram
\[\xymatrix{\widetilde{ S\times S}\ar[r]^r \ar[d]^g & S\times S \ar[d] &\ar[l]_{\;\;\;\; f} S\\
S^{[2]} \ar[r]^{r'}&\Sym^2 S,}
\]
where $r$ is the blowup of $S\times S$ along $f(S)$ and $f$ is  the diagonal embedding $S\ra S\times S$.
Let $S^{[2]}$ be the minimal resolution of $\Sym^2 S$ given by resolving the transversal $A_1$-singularity along
the diagonal and denote by $E$ the exceptional divisor of $r$. Let $\Delta$ denote $g_*(E)_{\red}$. The analogue
of the following proposition is well known for singular cohomology (see \cite[Proposition 6]{BeauKahl} ).

\begin{Prop} \label{symcoh}The diagonal embedding $f$ induces an isomorphism of Galois modules
\[ H^2_{\et}(S^{[2]},\Q_{\ell})\cong H^2_{\et}(S,\Q_{\ell})\oplus \Q_{\ell} [\Delta] .\]
\end{Prop}
\begin{proof} Since $S$ is a $K3$ surface we have that $H^1_{\et}(S,\Q_{\ell})=0$. Hence the
K\"unneth decomposition implies
\[H^2_{\et}(S\times S,\Q_{\ell})\cong H^0_{\et}(S,\Q_\ell) \otimes H^2_{\et}(S,\Q_{\ell})\oplus
H^2_{\et}(S,\Q_\ell)\otimes H^0_{\et}(S,\Q_{\ell}).\]
Let $\iota$ be the involution of $S\times S$ sending $(x,y)$ to $(y,x)$. Clearly $\iota^*$ acts
on $H^2(S\times S,\Q_{\ell})$ as (with respect to the K\"unneth decomposition)
\[ \left( \begin{array}{cc} 0&I_{22} \\I_{22}&0 \end{array}\right).\]
Hence $\iota^*$ has 22 eigenvalues 1 and 22 eigenvalues $-1$. Since $\iota^*$ fixes the 22-dimensional
module $H^2_{\et}(S,\Q_{\ell})$ (embedded via $f$) we obtain that
\[ H^2_{\et}(S\times S,\Q_{\ell})^{\iota^*}\cong H^2_{\et}(S,\Q_{\ell}).\]
Let $\tau$ be the involution on $\widetilde{S\times S}$ associated to $g$. One easily sees that
\[H^2_{\et}(S^{[2]},\Q_{\ell})\cong  H^2_{\et}(\widetilde{S\times S},\Q_{\ell})^{\tau^*}\cong
H^2_{\et}(S,\Q_\ell)\oplus \Q_{\ell}[\Delta]. \]
\end{proof}

\begin{Cor} \label{fourfoldcoh}Let $S/K$ be a K3 surface. Let $X/K$ be a cubic fourfold.
Assume that $F(X)\cong S^{[2]}$. Then
\[ H^4_{\et}(X,\Q_{\ell}) \cong H^2_{\et}(S,\Q_{\ell})(1) \oplus \Q_{\ell}[\Delta](1) \]
as Galois-modules.
\end{Cor}
In Section~\ref{sec_lseries} we prove that there exist $K3$ surfaces $S$ satisfying the assumptions of
this corollary.
Using this corollary and the Hodge conjecture (which is proven for cubic fourfolds), one obtains that if
$F(X)\cong S^{[2]}$ then $X$ is a special cubic fourfold.
Using lattice theory one can show that if $X$ is special of discriminant $d$ then $S$ admits a polarization of
degree $d$. (See \cite{HasArt}.)

\section{Moduli of cubic fourfolds}\label{sec_moduli}

The coarse moduli space $\cC$ of cubic fourfolds is the GIT-quotient
$\cC = V// \SL(6,\C)$, where $V \subset \Ps(H^0(\Ps^5, \str_{\Ps^5}))$
is the Zariski-open subset of cubic equations whose zero-locus is smooth.
Recall that all points in $V$ are properly stable in the sense of GIT. The variety $\cC$ is a
quasi-projective variety of dimension $20$. We define $\cC_d$ as the set of all
special cubic fourfolds $X$ that contain a surface $T$ such that $d(X,T)=d$.

\begin{Prop}\label{prpExistence}
Let $d>6, d \equiv 0,2 \bmod 6$. Then $\cC_d$ is a non-empty irreducible divisor in the variety $\cC$.
\end{Prop}
\begin{proof}
This is \cite[Theorem 1.0.1]{HasArt}.
\end{proof}

We say that a $K3$ surface $S$ is of \emph{degree $d$}
if $S$ has a polarization of degree $d$ (here $d$ is necessarily even).

\begin{Thm}\label{prpExistence2}
Assume that $d$ equals $2(n^2+n+1)$ where $n$ is an integer at least $2$. Then there exists an open set
$\cU_d$ of $\cC_d$ such that for every $X \in \cU_d$ there exists a $K3$ surface $S$ of degree $d$ such that
$F(X)\cong S^{[2]}$.
\end{Thm}
\begin{proof}
This is \cite[Theorem 1.0.3]{HasArt}.
\end{proof}

\section{Descent}\label{sec_lseries}
Fix a number field $K$. There are many examples of fourfolds $X$ such that $F(X)_{\C}\cong S^{[2]}$
where $S/\C$ is a  $K3$ surface and $S^{[2]}$ is the minimal desingularization
of $\Sym^2 S$ (see Proposition~\ref{prpExistence}). Suppose $S$ is defined over $K$.
In this section we prove that then both $X$ and the isomorphism $S^{[2]}\ra F(X)$ descend to $K$, provided
that $S^{[2]}(K)\neq \emptyset$ and the associated line bundle to $S^{[2]}_\C\cong F(X)_{\C}$ descends.

\begin{Def}
Let $T$ be a Noetherian scheme. A {\em family of quotients} of $\str_T^r$ parameterized by $T$ is a coherent sheaf $\cF$ flat over $T$, together with a surjective $\str_T$-linear homomorphism $q: \str_T^r\ra \cF$. Two pairs $(\cF,q)$ and $(\cF',q')$ are called isomorphic if $\ker(q)=\ker(q')$.   Consider the functor from (Noetherian) schemes to sets
\[ T \mapsto \left\{\begin{array}{c} \mbox{Isomorphism classes of locally free  quotients of }
\str_T^{n+1} \stackrel{q}{\ra} \cF,\\ \mbox{ parameterized by $T$, with } \rank \cF=k+1.\end{array}\right\}.\]
This functor is representable (\cite{Ni}) by a scheme which we denote by $\Gr(k,n)$. We  call $\Gr(k,n)$  the
Grassmannian of $k$-dimensional subspaces in $\Ps^n$.

A {\em vector bundle} on $X/T$ is a locally free sheaf on $X$, flat over $T$.
\end{Def}

Suppose $L$ is a field. Then one easily shows that the variety $\Gr(r,n)_L$ parameterizes dimension $r$ linear
subspaces of $\Ps^n_L$.

The following classical result describes all morphisms to Grassmann varieties:

\begin{Prop}\label{grasprp}Let $R$ be a ring with unit. Let $Y$ be a projective scheme over $\spec R$.
Giving a morphism $\phi: Y\ra \Gr(k,n)_R$  is equivalent to giving a rank $k+1$ vector bundle $V$ on $Y$
together with a surjective morphism of schemes $R^{n+1}\otimes \str_Y\ra V$ up to multiplication by an element in $R$. The composition of $\phi$ with
the Pl\"ucker embedding is given by the morphism associated to $\wedge^{k+1} R^{n+1} \otimes \str_Y\ra \wedge^{k+1} V$.
\end{Prop}
\begin{proof} See \cite[Proposition 2.1]{Ar} for the case that $R$ is a field.
The general case is a formal consequence of the fact that Quot-functors on Noetherian relative projective schemes
are representable. For this formal deduction we refer to \cite[Example 2.7]{Vi}, for a proof of the
representabiliy of the Quot-functor we refer to \cite{Ni}.
\end{proof}

\begin{Lem}\label{lemmaxdim} 
Let $L$ be a field. Fix positive integers $k,n$. Suppose $2k<n-1$.
Let $m= \left( \begin{array}{c} n+1\\k+1 \end{array} \right)-1$. Suppose $P$ is a plane contained in the image of the Pl\"ucker embedding of $\Gr(k,n)_K$ in $\Ps^m_K$. Then $\dim P\leq n-k$.

Moreover if $\dim P=n-k$ then $P$ corresponds to the family of $k$-dimensional planes containing a fixed  $k-1$-di\-men\-sio\-nal plane. Conversely, every $k-1$-di\-men\-sio\-nal plane gives such a family.
\end{Lem}

\begin{proof}
 The first claim is an easy excercise in Pl\"ucker relations. 

The second claim is an exercise in decomposable tensors: we illustrate this for the case $k=1$. In this case we need to show that a maximal linear subspace has dimension $n-1$ in $\Gr(1,n)$, and that this subspace correspond to lines throug a fixed point.

First let $v_1$ and $v_2$ be two decomposable $2$-tensors.
If the corresponding lines $l_1$ and $l_2$ in $\Ps^n_L$ are skew, then we 
can assume $v_1=e_0 \wedge e_1$ and $v_2=e_2 \wedge e_3$.
But then $v_1 \wedge v_2 \neq 0$ and the Pl\"ucker relations are, therefore, not satisfied, i.e., the line
spanned by $v_1$ and $v_2$ is not in the Grassmannian. For this to happen, we must have 
$v_1=e_0 \wedge e_1$ and $v_2=e_0 \wedge e_2$ (after choosing a suitable basis). 
Assume we have a third
decomposable $2$-tensor $v_3$ such that every tensor in the plane spanned by $v_1,v_2$ and $v_3$ is 
decomposable. Then the line $l_3$ corresponding to $v_3$ must meet $l_1$ and $l_2$. This can happen in two ways.
If $l_3$ is in the plane spanned by $l_1$ and $l_2$, then we can asusme $v_3=e_1 \wedge e_2$. This gives
us a plane in the Grassmannian, namely the plane parametrizing all lines in a given $\Ps^2_L$. If 
$n=3$ this is indeed a maximal dimensional subspace of $\Gr(1,3)_L$, but then $2k=n-1$. Again using the 
Pl\"ucker relations one sees that the plane spanned by $v_1=e_0 \wedge e_1, v_2=e_0 \wedge e_2$ and 
$v_3=e_1 \wedge e_2$ is not contained in any higher dimensional linear subspace contained in the
Grassmannian. The second possibility 
for $l_3$ is that it meets the point of intersection of $l_1$ and $l_2$, in which case
$v_3=e_0 \wedge e_3$, again in a suitable basis. We can now continue this argument to conclude that 
any $(n-1)$-dimensional
linear subspace of the Grassmannian is spanned by tensors of the form $e_0 \wedge e_1, \ldots ,
e_0 \wedge e_n$, i.e, all lines containing the point given by $e_0$.
\end{proof}

\begin{Lem}\label{grasaut} Let $L$ be a field. Fix non-negative integers $k, n$. Suppose that $2k< n-1$. Then $\Aut(\Gr(k,n)_L)=\Aut(\Ps^n_L)$.\end{Lem}
\begin{proof}
Let $m= \left( \begin{array}{c} n+1\\k+1 \end{array} \right)-1$ and let $\varphi: \Gr(k,n)_K\ra \Ps^m_K$ be the 
Pl\"ucker embedding.


We shall prove the lemma by induction on $k$. The statement is clear for $k=0$. {}From Lemma~\ref{lemmaxdim} it follows that 
an automorphism of $\Gr(k,n)_L$ induces an automorphism of $\Gr(k-1,n)_L$ by considering its action on maximal dimensional subspaces. 
By induction, this automorphism is induced from an automorphism of $\Ps^n_L$.
Now suppose that we have an automorphism $\sigma$ of $\Gr(k,n)_L$ inducing the identity on $\Gr(k-1,n)_L$. 
We claim that $\sigma$ is the identity on $\Gr(k,n)_L$. Indeed, consider a subspace $\Ps^k_L$. For any subspace $\Ps^{k-1}_L$ we know that $\sigma(\Ps^k_L)$ also contains $\Ps^{k-1}_L$. Hence $\sigma(\Ps^k_L)=\Ps^k_L$.
\end{proof}

\begin{Rem} The Grassmannian $\Gr(1,3)_K$ contains two 
families of planes (they are usually called $\alpha$-planes and $\beta$-planes). There
is an automorphism of the Grassmannian $\Gr(1,3)_K$ (which is a quadric in $\Ps^5_K$), which 
interchanges the two types of planes and is, therefore, not induced by a linear map of $\Ps^3_K$.
\end{Rem}

Fix a $K3$ surface $S/K$. Assume that there exists a number field $L\supset K$ and a cubic fourfold $X/L$ such
that $F(X)\cong S^{[2]}_{L}$, where $F(X)\subset \Gr(1,5)$ denotes the Fano variety of lines of $X$.

For the rest of the section fix a rank 2 vector bundle $V$ on $S^{[2]}_L$ yielding an
isomorphism $S^{[2]}_L\cong F(X)\subset \Gr(1,5)$. We proceed now by proving that $V$ descends to $S^{[2]}$ under the assumption that $\det(V)$ descends to $S^{[2]}$.




{}From here on we will assume that $\det(V)$ on $S^{[2]}_L$ is isomorphic to the pullback of a line bundle
on $S^{[2]}$. We express this by saying that {\em $\det(V)$ is defined over $K$}. This is equivalent to stating that the associated line bundle (cf. Section~\ref{intro}) on $S$ descends to $K$. One can show that $\det(V)$ is a combination of the associated line bundle and $\delta$, the square root of the exceptional divisor $\Delta$. So we need to prove that $\delta$ is defined over $K$. For any element $\sigma \in \Gal(\overline{K}/K)$ we have $\delta\otimes \delta=\Delta=\delta^\sigma\otimes \delta^\sigma$. Since $\Pic(S^{[2]})_{\tor}$ is trivial this implies that $\delta^\sigma \cong \delta$. Since $S^{[2]}(K)$ is non-trivial we can `rigidify' $\delta$. Using descent theory as in \cite[Theorem 2.5]{Klei} one can show that $\delta$ is defined over $K$.

To assume that $\det(V)$ is defined over $K$  is not a severe restriction: one can show that for all examples
produced by Hassett \cite{HasArt}  $\det(V)$ is a linear combination of the polarizing class of $S$ and the square root $\delta$. 
The exceptional divisor is defined over $K$. 
%
 In general it is not hard to check whether the polarizing class is defined over $K$.




\begin{Dia}\label{basechngdia} Let $R:=L\otimes_K L$ and $R':=L\otimes_K L\otimes_K L$.
Denote by $h_i: \spec R \ra \spec L$ and by $h_{j,k}: \spec R' \ra \spec R$   the morphism induced
by the inclusion on the $i$-th and on the $j$-th and $k$-th factor, respectively. Similarly, denote by
$s_n: \spec R' \ra \spec L$ the morphism induced by the inclusion of the $n$-th factor.
The morphism $p_i$ is the morphism $S^{[2]}_R\ra S^{[2]}_L$ induced by the base change $h_{i}$,
the morphism $p_{j,k}$ is defined analogously. Set $n=j$ if $i=1$, or $n=k$ if $i=2$. Then we obtain the
following commutative diagram:
\[\xymatrix{ S^{[2]}_{R'}  \ar[d] \ar[r]_{p_{j,k}} \ar @/^/ [rr]^{q_n}& S^{[2]}_R  \ar[d] \ar[r]_{p_i} & S^{[2]}_L \ar[d]
\ar[r]_{p}& S^{[2]} \ar[d]
\\
\spec R'  \ar[r]^{h_{j,k}} \ar @/_/[u]_{g_{R'}} \ar @/_/ [rr]_{s_n} & \spec R  \ar[r]^{h_{i}}
\ar @/_/[u]_{g_{R}} &\spec L \ar[r]^{h}\ar @/_/[u]_{g_{L}} &\spec K.\ar @/_/[u]_{g}}\]
\end{Dia}

\begin{Rem} Suppose for the moment that $L/K$ is Galois. Write $L:=K[t]/(f)$, for some polynomial $f\in K[t]$
of degree $d$. Then
\[ R=L\otimes_K L = K[t]/(f) \otimes_K L = L[t]/(f)\cong L^{d},\]
where the isomorphism follows from the Chinese remainder theorem.
In particular, $\spec R$ consists of $d$ closed points.

One can show that we can choose the isomorphism $L\otimes_K L\ra L^d$ such that
\[ w\otimes 1\mapsto (w,w,\dots,w) \mbox{ and } 1\otimes w\mapsto (w,g_1(w),g_2(w),\dots, g_{d-1}(w)),\]
where $\Gal(L/K)=\{1,g_1,\dots, g_{d-1}\}$.

Similarly, one can show that $S^{[2]}_R$ is the disjoint union of $d$ copies of $S^{[2]}_L$.
\end{Rem}

\begin{Rem}Suppose we sheafify the functor $B$ from $K$-schemes to sets given by
\[ T\mapsto \{ \mbox{Isomorphism classes of locally free sheaves of rank } r \mbox{ on }X_T \}/\Pic(T)\]
in the \'etale, fpqc or fppf topology. Denote the sheafified functor by $B'$. The following lemma
shows that $V\in B'(\spec K)$. (In this context, our aim is to prove that $V\in B(\spec K)$.)
\end{Rem}




\begin{Rem} Let $f:S^{[2]}\ra \spec K$ be the structure map. We say that $f_*\str_{S^{[2]}}\cong \str_{\spec K}$
holds universally, if there is an isomorphism $\psi : f_*\str_{S^{[2]}}\ra \str_{\spec K}$ such that for any
morphism of schemes $h:T\ra \spec K$ we have that the base-changed morphism
of $\str_T$-modules $\psi_T : f_{T*} \str_{S_T^{[2]}} \ra \str_T$ is an isomorphism.

In \cite[Answer to exercise 3.11]{Klei} it is shown that if $X/B$ is a proper and flat family, with connected and reduced fibers,
then $f_*\str_X \cong \str_B$  holds universally. Since we are working with a projective variety over a
field, we are in this situation.
\end{Rem}


%
%

\begin{Prop}\label{descentPrp} Suppose $S^{[2]}(K)\neq \emptyset$ and $\det V$ descend to $S^{[2]}$. Then the vector bundle $V$ on $S^{[2]}_L$
descends to $S^{[2]}$.\end{Prop}

\begin{proof}
Assume the notation from Diagram~\ref{basechngdia}. Without loss of generality we may assume that $L/K$ is Galois. This implies that $\spec R\ra \spec L$ is a finite union of points.

We start by constructing an isomorphism $w:p_1^*V\ra p_2^*V$. Consider the embedding $\varphi_i$ associated with $\wedge^2 p_i^*V$. Since $\wedge^2 p_1^*V\cong \wedge^2 p_2^*V$ we obtain that the image of  $\varphi_1$ and $\varphi_2$ differ by an automorphism of $\Ps^{14}_R/\spec(R)$ fixing the Pl\"ucker embedding of $\Gr(1,5)_R$. Since $\Gr(1,5)_R=\sqcup \Gr(1,5)_L$, we obtain from  Lemma~\ref{grasaut} that $\Aut(\Gr(1,5)_R/\spec(R)\cong \Aut(\Ps^5_R/\spec R)=\PGL_6(R)$. In particular there is a commutative diagram
\[ \xymatrix{ \str^6_{S^{[2]}_{R'}} \ar[r]\ar[d]^{w_1} &p_1^*V \ar[r]\ar[d]^{w}& 0\\\str^6_{S^{[2]}_{R'}} \ar[r]&p_2^*V \ar[r]& 0
}\]
where the vertical arrows are isomorphisms.

By descent theory (e.g., \cite[Theorem I.2.23]{Mil})
it suffices to show that the automorphism 
\[\tau:= (p_{13}^*w)^{-1}\circ p_{23}^*w\circ p_{12}^*w \]
 of $q_1^*V $ is trivial. 

Consider the quotient
\[ s: \str_{S^{[2]}_{R'}}^6\longrightarrow V_{R'}\]
induced by  the embedding $S^{[2]}_{R'}\ra \Gr(1,5)_{R'}$. Then $\tau \circ s$ gives another embedding $S^{[2]}_{R'}\ra \Gr(1,5)_{R'}$. As above we can compose these two embedding with the Pl\"ucker embedding and we get two embeddings into $\Ps^{14}_{R'}$ that differ by an automorphism of $\Gr(1,5)_{R'}/\spec R'$. Each such automorphism is coming from an automorphism of $\Ps^5_{R'}$.
 Equivalently, we have a diagram
\[ \xymatrix {\str_{S^{[2]}_{R'}}^6 \ar[r] \ar[d]_{\tau_1} & V  \\ \str_{S^{[2]}_{R'}}^6\ar[ur]}\]
where $\tau_1$ is an element of $\GL_6(R')\subset H^0(\GL_6(\str_{S^{[2]}_{R'}}))$.


Fix a point in $S^{[2]}(K)$ and let $g: \spec K\hookrightarrow S^{[2]}$ be the corresponding morphism. Then by the definition of fibre products  and the fact that $g_{R'}$ corresponds to a $K$-rational point we obtain that $g_{R'}^*\tau_1$ is trivial.
Now, 
\[\tau_1 \in H^0(\GL_6(\str_{S_{R'}^{[2]}})) = H^0(\GL_6(f_*\str_{S_{R'}^{[2]}}))\stackrel{\psi'_{R'}}{\ra} H^0(\GL_6(\str_{\spec R'}))
\mbox{ (as groups)}.\]
Since $\psi_{R'}: f_*\str_{S_{R'}^{[2]}}\cong \str_{\spec R'}$, we have that the induced map $\psi'_{R'}$  is an
isomorphism of groups. Since $g_{R'}^*\tau_1$ is trivial we obtain that $\psi'_{R'}(\tau_1)$ is also trivial.
This implies that $\tau$ is trivial, hence   $V$ descends to $S^{[2]}$.
\end{proof}

\begin{Thm}\label{mainthm}
Let $S/K$ be a $K3$ surface such that $S^{[2]}(K)$ is non-empty and such that over some finite
extension $L/K$
we have $S^{[2]}_L\cong F(X)$ with $X/L$ a cubic fourfold. Suppose the associated line bundle (cf. Section~\ref{intro}) descends to $K$.
Then the cubic fourfold $X$ can be defined over $K$ in such a way that $S^{[2]}\cong F(X)$.\end{Thm}
\begin{proof}
It follows from Proposition~\ref{descentPrp} that $V$ is vector bundle on $S^{[2]}$.
Consider the evaluation map
\[ \omega: H^0(V)\otimes \str_{S^{[2]}} \ra V.\]
Since forming cohomology commutes with flat base change, and $\omega_L$ is surjective, it follows
that $\omega$ is surjective. This implies that the vector bundle $V$ is globally generated.

By Proposition~\ref{grasprp} we obtain that the associated map $\psi:S^{[2]}\hookrightarrow \Gr(1,5)$
is defined over $K$.
Let $U\subset \Gr(1,5)\times \Ps^5$ be the universal line and $p$ and $q$ the projections.
If $X'/K$ is the image of $p^{-1}\psi(S^{[2]})$ under $q$, then, by construction, we have $X'_L\cong X$.
Hence $X'$ is a model of $X$, so $X$ has a model over $K$ and $\psi:S^2\stackrel{\sim}{\ra} F(X')$.
\end{proof}

\begin{Cor}\label{isocor} Fix $n\geq 2$ and let $d=2(n^2+n+1)$.
There is a non-empty Zariski open set $\cV_d$ of the moduli space $\cF_d$ of polarized $K3$ surfaces
of degree $2d$ such that if
\begin{enumerate}
\item $K$ is a subfield of $\C$;
\item $S/K$ a $K3$ surface such that $S_\C$ corresponds to a point
in $\cV_d$;
\item the set $S^{[2]}(K)$ is not empty;
\item all ample line bundles of degree $d$ on $S$ descend to $K$;
\end{enumerate}
then there exists a cubic fourfold $X/K$
such that
$F(X) \cong  S^{[2]}$.
\end{Cor}
\begin{proof} By \cite[Theorem 1.0.3]{HasArt} there is a set $\cV_d$ such that for every $K3$ surface $S/\C$
corresponding to a point in $\cV_d$ there exists  a cubic fourfold $X/\C$ such that $S^{[2]}\cong F(X)$.

Assume $S$ can be defined over $K$. Consider the subscheme $T$ of the Hilbert scheme of $\Gr(1,5)$
corresponding to subvarieties of $\Gr(1,5)$ that are geometrically isomorphic to  $S^{[2]}$, and the isomorphism can be given by an element of $\Aut(\Ps^5)(\overline{\Q})\subset \Aut(\Gr(1,5))(\overline{\Q})$.
It easy to see that $T$ is defined over $K$. By assumption $T$ is not empty, hence  $T(\overline{K})\neq \emptyset$.
In particular, there exists a finite extension $L/K$ for which $T(L)\neq \emptyset$. Hence $S^{[2]}_L$ can be
embedded in $\Gr(1,5)_L$. Since $qp^{-1}(S^{[2]}_L)$ is a cubic fourfold, we obtain that $S^{[2]}_L\cong F(X)$
for some cubic fourfold $X/L$. From \cite[Theorem 1.0.3]{HasArt} it follows that the associated line bundle is an ample line bundle of degree $d$. Since all ample line bundles of degree $d$ are defined over $K$, we have that the associated line bundle to $F(X)\cong S^{[2]}$ is defined over $K$. Now apply Theorem~\ref{mainthm}.
\end{proof}

We can now relate the zeta function of the cubic fourfold $X$ to the zeta function of the $K3$ surface $S$.

\begin{Cor}\label{lsercor} Let $S,X,K$ be as before. Then we have the following equality of zeta functions
\[ Z_K(X,s)\circeq Z_K(S,s-1) \zeta_K(s)\zeta_K(s-2)\zeta_K(s-4).\]
Here $\circeq$ means equality upto, possibly, finitely many local Euler factors.
\end{Cor}

\begin{proof}
First recall that $H^1(X,\Z)=H^3(X,\Z)=0$ and that $H^2(X,\Z)=\Z h$ where $h$ is the class of the hyperplane section.
Moreover $H^1(S,\Z)=0$.  It follows from Corollary~\ref{fourfoldcoh} that
\[H^4_{\et}(X,\Q_{\ell})\cong H^2_{\et}(S,\Q_{\ell})(1) \oplus \Q_{\ell} [\Delta](1).\]
Hence we have that
\[\frac{Z_K(X,s)}{\zeta_K(s)\zeta_K(s-1)\zeta_K(s-3)\zeta_K(s-4)}\circeq
\frac{ Z_K(S,s-1)\zeta_K(s-2)}{\zeta_K(s-1)\zeta_K(s-3)}\]
which gives the claim.
\end{proof}

If the surface $S$ is singular, we can determine the zeta function explicitly.
Recall that a $K3$ surface $S/K$ is called \emph{singular} if its geometric Picard number $\rho(S_{\overline{K}})$
equals $20$. Singular $K3$ surfaces were classified
by Shioda and Inose \cite{SI}. The isomorphism classes of these surfaces are in $1$-to-$1$ correspondence with
$\SL(2,\Z)$-isomorphism classes of even binary quadratic forms. This correspondence is given by associating
to a singular $K3$ surface $S$ its transcendental lattice $T(S)$. The discriminant of the surface $S$ is
defined as the discriminant of the lattice $T(S)$. Shioda and Inose showed that every
singular $K3$ surface $S$ is either the
Kummer surface of a product $E \times F$, where $E$ and $F$ are isogenous elliptic curves with complex
multiplication (CM), or $S$ has a (rational) double cover by such a surface. The elliptic curve $E$ has
CM in the field $\Q(\sqrt{-d_0})$ where $d_0$ is the discriminant of the surface $S$. The construction of
Shioda and Inose allows one to have control over the field of definition $K$.

Using this classification, Shioda and Inose have computed the zeta function of singular $K3$ surfaces
(\cite[Theorem 6]{SI}). If the field $K$ is chosen big enough, namely such
that the N\'eron-Severi group can be generated by elements
defined over $K$ (this means essentially that the points of order $2$ of the curves $E$ and $F$ are
defined over $K$), then the zeta function of $S$ is a product whose factors are either (twisted)
Dedekind zeta functions or the $L$-function of a suitable Hecke Gr\"ossencharakter. For a recent account of
Hecke Gr\"ossencharakters and their associated Hecke eigenforms, in particular in connection with $K3$
surfaces, we refer the reader to \cite{Sch}. It is a consequence of our descent Theorem~\ref{mainthm} that there
is a cubic fourfold $X$ defined over $K$ with $S^{[2]}\cong F(X)$. As a consequence we can compute the
zeta function of the variety $X$.

We recall that the singular K3 surfaces form a countable set in the moduli space $\cF_{d}$ which is dense in the $\C$-topology. This can be proved in the same way as the density in the period domain of all $K3$ surfaces (see the proof of  \cite[Corollary VIII.8.5]{BHPV} ).
This shows, using Theorem~\ref{prpExistence2}, that for every value of $d=2(n^2+n+1), n \in\Z_{>1}$ 
there exists a countable number of points in a non-empty
Zariski open subset of $\cC_d$ such that the corresponding fourfolds have a model over an explicitly known
number field $K$ and where the $L$-function can be computed in terms of Hecke Gr\"ossencharakters.

\begin{Cor} \label{cormodA}Assume $S$ satisfies the assumptions of Corollary~\ref{isocor}. In particular, there
exists a cubic fourfold $X/K$ with $S^{[2]}\cong F(X)$. Moreover assume that
$\rho(S)=\rho(S_{\overline{K}})=20$, that the square root of the negative of the discriminant of the transcendental lattice of $S$
is contained in $K$ and that the N\'eron-Severi group of $S$ can be generated by divisors defined over $K$.
Then there exists a Gr\"ossencharakter $\psi$ such that
\[ Z_K(X,s)\circeq \zeta_K(s)\zeta_K(s-1)\zeta_K(s-2)^{21}\zeta_K(s-3)
\zeta_K(s-4)L(\psi^2,s-1)L(\overline{\psi}^2,s-1).\]
\end{Cor}
\begin{proof} {}From \cite[Theorem 6]{SI} it follows that there exists a Gr\"ossencharakter $\psi$ such that
\[ L(S,s)\circeq \zeta_K(s)\zeta_K(s-1)^{20}\zeta_K(s-2)L(\psi^2,s)L(\overline{\psi}^2,s).\]
Combine Corollaries~\ref{isocor} and~\ref{lsercor} to conclude the proof.
\end{proof}

If a singular $K3$ surface $S$ can be defined over $\Q$, then the transcendental lattice $T(S)$ is a
$2$-dimensional Galois module. Livn\'e \cite[Example 1.6]{LivMot} has shown that in this case the $L$-function defined by $T(S)$
is that of a Hecke eigenform of weight $3$. Assume that $S^{[2]}\cong F(X)$ for some cubic fourfold.
Then we can, by our descent theorem, assume that $X$ is defined over $\Q$. The sublattice of $H^4(X,\Z)$
generated by algebraic cycles has rank $21$. We denote its orthogonal complement by $T(X)$ and refer to this as
the {\em transcendental lattice} of $X$.

\begin{Cor}\label{cortranscendentallattice} Assume that $S$ is a singular $K3$ surface defined over $\Q$,
and that it
satisfies the assumptions of Corollary~\ref{isocor}, in particular $S^{[2]}\cong F(X)$ for some cubic
fourfold $X$ defined over $\Q$. Then  there exists a weight $3$ Hecke eigenform $f$ such that
\[ L(T(X),s)\circeq L(f,s-1).\]
\end{Cor}
\begin{proof}
 Proposition~\ref{prpGaloisisomorphism} implies that 
$T(X) $ and $T(S)(1)$ are isomorphic.
\end{proof}

\section{The Fermat cubic}\label{sec_fermat} 

We shall now present an explicit example, showing that the weight $3$ form associated to the
Fermat cubic is the level $27$ newform. 

Our example will be a special case of the Pfaffian construction due to Beauville and Donagi \cite{BeDo}.
For this we shall consider a cubic fourfold $X$
containing two planes $P$ and $Q$, which we can assume to be $P=\{u=v=w=0\}$ and $Q=\{x=y=z=0\}$.

The defining polynomial for $X$ can be written in the form $F-G$ where $F,G$ are
bihomogeneous in $\C[u,v,w;x,y,z]$ of bidegrees $(2,1)$ and $(1,2)$ respectively. Such a fourfold is
always rational. Indeed, a general line joining the two planes $P$ and $Q$ will meet $X$ in a third point.
In this way we obtain a birational map $\pi_1\times \pi_2:  X\dashrightarrow \Ps^2\times \Ps^2$.
The lines contained in $X$ correspond to the complete intersection $S$ given by $F=G=0$. The
surface $S$ is a complete intersection of bidegree $(2,1), (1,2)$ and hence a $K3$-surface.
Via the Segre embedding $\Ps^2 \times \Ps^2\ra \Ps^8$ it is embedded into
$\Ps^8$, where it is of degree $14$. This corresponds to the case $n=2$ in Theorem~\ref{prpExistence2}.

The inverse birational map is given by the linear system $|L|$ of bidegree $(2,2)$ forms on
$\Ps^2\times \Ps^2$ vanishing along $S$. In this case, as was shown by Hassett (see \cite[page 41]{Has}),
the map $S^{[2]}\to F(X)$ which was described in \cite[Proposition 5]{BeDo} is no longer an isomorphism.
However, there is a birational map
$S^{[2]}\dashrightarrow F(X)$, which will be good enough for our purposes:
Take a point $s+t\in S^{[2]}$ and consider the lines $l_i$ spanned by $\pi_i(s)$ and
$\pi_i(t)$ for $i=1,2$.
If $s$ and $t$ are sufficiently general, then $l_1\times l_2\subset \Ps^2\times \Ps^2$
intersects $S$ in precisely five points. Two of these points are $s$ and $t$, call the other three
points $u,v,w$. Let $C$ be the $(1,1)$ curve on $l_1\times l_2$ passing through $u,v,w$.
Since the three points $u,v,w$ lie on $S$, the linear series $|L|$ restricted to $C$ has
degree $1$, and hence
$C$ is mapped to a line in $\Ps^5$.
One can show that the general line on $X$ is obtained in this way, giving the desired birational map.

In order to present an explicit example, consider the surface $S$, given by the complete intersection of
\[ F=xu^2+yv^2+zw^2=0 \mbox{ and } G=x^2u+y^2v+z^2w=0.\] 
This surface was considered in \cite{DoKo} as an example for a $K3$-surface in characteristic $2$ such that
the automorphism group has a quotient that is not isomorphic to a subgroup of the 
Mathieu group $M_{23}$, proving that Mukai's classification of finite groups acting symplectically on 
complex $K3$ surfaces does not extend to characteristic $2$.

The surface $S$ contains many lines. An easy calculation, using the intersection pairing,
shows that $\rho(S)=20$. Moreover, the discriminant of the transcendental lattice $T(S)$ 
equals $3$ up to a square.
One can also show that $3$ is the only prime of bad reduction.
In particular, $L(S,s)\circeq L(f,s-1)$, with $f$ a weight $3$ Hecke newform of level $3^b$, for some $b$. There are precisely three such forms, namely the newform of level $27$ and two twists of this newform, by the standard cubic character and the square of the standard cubic character. The Fourier coefficients at $7$ of these three forms are pairwise distinct.

We will determine $b$ and $f$: Let $p$ be a prime of good reduction and let $\overline{S}$ be the reduction 
of $S$ modulo $p$. Consider the action of Frobenius $\Fr$ on the cohomology of $\overline{S}$. On $H^4_{\et}(S,\Q_\ell)$ Frobenius acts as multiplication by $p^2$, on $H^0_{\et}(S,\Q_\ell)$ as multiplication by $1$. 
It is well-known to the experts  that the trace of Frobenius on \[NS(\overline{S})\otimes \Q_\ell \subset H^2_{\et}(\overline{S},\Q_{\ell})\] is 0 modulo $p$. From the Lefschetz trace formula it follows that
\[ \Tr \Fr^* \mid T(\overline{S})\equiv \#\overline{S}(\F_p) -1 \bmod p.\]

For primes $p\equiv 2 \bmod 3$, one easily shows, using the form of the equations, that 
this trace is $0 \bmod p$. For primes $p\equiv 1 \bmod 3$ a straightforward calculation by computer gives:
\[\begin{array}{c|ccccc} p&7 & 13 & 19 & 31 &37\\\hline
\#\overline{S}(\F_p)&177& 429 &753 &1536 & 2157 \\
(\Tr \Fr^*\mid T(\overline{S}))  \bmod p&1 & 12 &11 & 16 & 10\\\end{array}\]
Since the form $f$ corresponds to the Galois representation of the transcendental lattice, 
the Fourier coefficients of $f$ coincide with $\Tr \Fr^*\mid T(\overline{S})$.
Using the tables from \cite{Sch} one concludes that $f$ is the unique newform of level $27$ up to a twist by a cubic character, unramified outside 3. By class field theory there are three such characters, say $\chi, \chi^2$ and the trivial character. One easily sees that the Fourier coefficients of $f$, $f\otimes \chi$ and $f\otimes \chi^2$ at 7 are different. From this we conclude that $f$ is the correct form. 

The corresponding fourfold $X$ can easily be determined explicitly. The linear system of
forms of bidegree $(2,2)$ which vanish on $S$ defines the map
\[
(xF:yF:zF:uG:vG:wG) : \Ps^2 \times \Ps^2 \dashrightarrow \Ps^5
\]
and an easy calculation shows that the image satisfies the equation
\[ F-G = xu^2-ux^2+yv^2-vy^2+zw^2-z^2w\]
which is the equation of the associated fourfold $X$. The birational map
$\psi: S^{[2]}\dashrightarrow F(X)$
induces an isomorphism $T(S)\cong T(F(X))$ of Galois modules. (This follows from the fact that $\psi$ induces an isomorphism
$\psi^*:H^{2,0}(S) \oplus H^{0,2}(S) \ra H^{2,0}(F(X)) \oplus H^{0,2}(F(X))$.)
As in Corollary~\ref{cortranscendentallattice} it now follows that
\[
L(T(X),s) \circeq L(f,s-1)
\]
where $T(X)$ is the transcendental lattice of $X$ and $f$ is the weight $3$ Hecke eigenform of level $27$
from above.

Since the equation of $X$ is the sum of three terms which only depend on two variables each, it follows
that $X_{\Q(\sqrt{-3})}$ is isomorphic to the Fermat cubic $X'$
\[u^3+v^3+w^3+x^3+y^3+z^3.\]
This isomorphism does not descend to $\Q$. 
This can be shown as follows. The automorphism group of $X'$ is generated by permutaions of the coordinates and by multiplying the coordinates by a third root of unity (see \cite{Shioda}). This implies that $\#\Aut(X')(\Q)=\#S_6=720=6^2\cdot 20$. One can easily find many $\Q$-rational automorphisms of $X$. For example  the automorphisms $[x,u,y,v,z,w]\mapsto [-u,-x,y,v,z,w]$ and $[x,u,y,v,z,w]\mapsto [u-x,-x,y,v,z,w]$ generate a subgroup of order $6$. Considering the same automorphisms for the pairs $(y,v)$ and $(z,w)$ gives a subgroup of order $6^3$ in $\Aut(X)(\Q)$. In particular, $\#\Aut(X)(\Q)\neq \Aut(X')(\Q)$. This implies that $X$ and $X'$ are not isomorphic
and that the $K3$ surface associated to $X'$ is a twist of our surface $S$. 
However, we shall show that the associated newform to $X'$ is the same form $f$. 

The Galois representation on $T(X')$ is a twist of the Galois representation on $T(X)$ by the quadratic
character 
$\chi_{-3}$. In particular, the Hecke newform $f'$ associated to $X'$ has level $3^c$, so it is 
one of the three above mentioned forms. From $\#\overline{X'}(\F_7)=3690$ one deduces, as above, that $f'=f$.

The zeta function of the Fermat cubic was also conidered by Goto in \cite{Go}. His approach, however,
is entirely different, since he considered this over the field $\Q(\sqrt{-3})$ and he used Weil's classical approach using Jacobi sums.

\end{document}